\documentclass[11pt]{article}
\usepackage{amsmath,amssymb}
\usepackage{enumerate}
\usepackage{color}
\usepackage{amsthm}

\newcommand{\R}{\mathbb{R}}

\newcommand{\N}{\mathbb{N}}
\newcommand{\bea}{\begin{eqnarray}}
\newcommand{\eea}{\end{eqnarray}}
\newcommand{\la}{\lambda}

\def\a{\alpha}

\def\e{\varepsilon}

\def\s{\sigma}

\def\x{\otimes}

\def\diam{{\rm diam}}
\def\ric{{\rm ric}}
\def\1{\rm Id}
\def\sgn{\rm sgn}
\def\sup{{\rm sup}}
\def\vol{\rm vol}

\def\ci{\circ}
\def\cl{{\rm cl}}
\def\pr{{\rm pr}}

\def\V{\noindent}
\def\Hess{{\rm Hess}}

\def\csec{{\rm csec}}
\def\cdiam{{\rm cdiam}}
\def\sec{{\rm sec}}

\def\injrad{{\rm injrad}}
\def\Vol{{\rm Vol}}

\newcommand{\bean}{\begin{eqnarray*}}
\newcommand{\eean}{\end{eqnarray*}}

\newtheorem{Theorem}{Theorem}

\newtheorem{Lemma}{Lemma}

\newcommand{\ben}{\begin{enumerate}}
\newcommand{\een}{\end{enumerate}}
\newcommand{\bit}{\begin{itemize}}
\newcommand{\eit}{\end{itemize}}
\newcommand{\edoc}{\end{document}}

\usepackage{hyperref}

\usepackage{enumerate}

\title{Lorentzian Gromov-Hausdorff theory and finiteness results}

\begin{document}

\author{Olaf M\"uller\footnote{Institut f\"ur Mathematik, Humboldt-Universit\"at zu Berlin, Unter den Linden 6, D-10099 Berlin, \texttt{Email: mullerol@math.hu-berlin.de}}}

\date{\today}
\maketitle

\begin{abstract}
\V We construct a functor from a category of Lorentzian manifolds-with-boundary to a category of Riemannian manifolds-with-boundary. As a corollary, we obtain the first known Cheeger-Gromov type finiteness result in Lorentzian geometry. 
\end{abstract}

\section{Introduction and statement of the main results}

Considering the tremendous success of Gromov's theory of metric spaces in Riemannian geometry \cite{mB}[Sec.14.6], it appears worthwhile to look for a Lorentzian version of this concept, more so as it is desirable --- e.g. in the context of non-perturbative Quantum Gravity --- to find good covariant ways to topologize spaces of Lorentzian metrics and their moduli spaces. A central object in Gromov's theory is the set $M(r,s ,D)$ of all {\em $(r,s,D)$-Alexandrov metrics} (for $(r ,s,D) \in \R^3 $), which are locally compact complete (thus geodesic) length spaces with Hausdorff dimension  $r$, curvature $ \geq s$ and diameter $\leq D $. The set $M(\leq r,   s,D ):= \bigcup_{u \leq r} M(u, s, D)$, if equipped with the Gromov-Hausdorff metric (defined first by Edwards \cite{Edwards}, for historical details see \cite{aT}), has been shown to be compact by Gromov \cite{mG}. Then Perelman's stability theorem (cf \cite{Kapovitch}) (stating that each converging sequence in $M(r,s, D)$ is eventually in one homeomorphism class) implies the finiteness result of \cite{GPW}.

As on Minkowski space $\R^{1,n}$ there is no metric compatible with the topology and covariant w.r.t. the Poincar\'e group (\cite{Invitation}, p.6, Item 3) (essentially because the orbit of every point in $\partial J^+(p)$ under a group of boosts accumulates at $p$), there is no injective functor $N$ from globally hyperbolic (g.h.) spacetimes to metric spaces preserving the topology, i.e., with $\mathcal{F} = \mathcal {F } \ci N $ for the forgetful functor $\mathcal{F}$ to topological spaces, so we restrict our categories: For $n \in \N$, let $C^-_n$ be the category of $n$-dimensional {\em Cauchy slabs}, i.e. ${\rm Obj} (C_n^-)$ consists of the g.h. $n$-dimensional manifolds-with-boundary $X$ whose boundary are two disjoint smooth spacelike Cauchy surfaces, the morphisms being isometric, oriented, time-oriented diffeomorphisms, and the category $C_n^+$ is the one of complete connected $n$-dimensional $C^2$ oriented Riemannian manifolds-with-boundary and isometries (those preserve the boundaries). To set up geometric bounds, we now define four real functions on the set of objects $(X,G)$ of $C_n^+ \cup C_n^-$ (so $G$ is of unspecified signature): 

\begin{itemize}
	\item $k_1: (X,G) \mapsto  \csec(X,G) := \inf \{  \sec (A) \vert  x \in X , A \subset T_xX  {\rm  \ linear \ two-dimensional \ subspace }, \\  A^\perp {\rm \ spacelike}\} $, inf of cospacelike sectional curvature (i.e. timelike sectional curvature on $C^-$);
	\item $k_2: (X, G) \mapsto \sup_{x \in \partial X} \{ \vert \nabla^G \nu (x) \vert_G : x \in \partial X  \} $ norm of second fundamental form (here $\nu (x)$ is the outer normal vector at $x$, past on $\partial^- X$ and future on $\partial^+ X$ if $(X,G) \in {\rm Obj} (C_n^-)$); 
	\item $k_3: (X,G) \rightarrow \Vol^G (X)  = \int_X \vol^G$;
	\item $k_4: (X,G) \mapsto \cdiam ((X,G)) := \sup \{ G(\nu (p) , w (p)) \vert p \in \partial X, G(w(p), w(p)) <0, \exp_p (w) {\rm \ exists }  \} $. For $G$ Riemannian, $\cdiam (X,G) = \diam (X,G)$, for $G$ Lorentzian, $\cdiam (X,G) \geq \sup \{ \ell (c) \vert c : \R \rightarrow X {\rm \ causal } \}$ by compactness of $ J^+(p)$ for each $p \in X$ and the first variational formula.
\end{itemize}

For objects $(X,g)$ of $C_n^-$ we moreover define, with $J_g(p) :=  J_g^-(p) \cup J_g^+(p)$:

\begin{itemize}
\item ${\rm Jvol} (X) := \sup \{ \vol (J (p) ) \vert p \in X \}, {\rm Jvol } (\partial X) = \sup \{ \vol (J (p) \cap \partial X  ) \vert p \in X \} , $
\item ${\rm injrad}^\pm_g (x) :=  \sup \{  \sqrt{- g(w,w)} : w \in T_xX, w \ll 0,  \exp\vert_ {I^\pm(0) \cap I^\mp(w)} {\rm \ is \ a \  diffeomorphism } \}  , $
\item $ \Gamma (X) := \inf_{x \in X} \max \{ {\rm injrad}^+_g (x) , {\rm injrad}^-_g (x)   \} .$
\end{itemize}

Then, for each $a \in (\R \cup \{ \infty \})^4$ and each $b \in (\R \cup \{ \infty \})^7$, we define subcategories

$$C_n^\pm (a) := \{ (X,g) \in C_n \vert - \csec (X) \leq a_1 \land \vert \nabla \nu \vert_{\partial X} \leq e^{a_2} \land  ({\vol} (X))^{-1} \leq e^{a_3} \land \cdiam (X) \leq e^{a_4} \} , $$

$$C_n^- (b) := \{ (X,g) \in C_n^- (b_1, b_2, b_3, b_4) )  \vert  {\rm Jvol}^g (X) \leq e^{b_5} \land (\Gamma (X))^{-1} \leq e^{b_6} \land {\rm Jvol}^g ({\partial} X ) \leq e^{b_7}  \} , $$

\V of $C_n ^\pm$. To be able to reconstruct the Lorentzian metric from the data on the Riemannian side (see the second item of Th. \ref{Main1}), we need a slightly richer category $K_n$ consisting of objects $(X,G)$ of $C^+_n$ carrying in addition a Lipschitz function $X \times X \rightarrow \R^{3}$, the morphisms acting by pull-back on them. We will find an injective functor $F: C_n \rightarrow K_n$ mapping appropriate subsets $ C^-_n(b) $ into some $C_n^+ (a)$. Let $\s^g: X \times X \rightarrow \R$  the signed Lorentzian distance $\sigma^g : X \times X \rightarrow \R \cup \{ \infty \} \cup \{ - \infty \}$ defined by $\s^g(x,y) := \mp \sup \{ \ell(c)  \vert c: x \leadsto y {\rm \  causal \  curve} \}$ for $x \in J^\pm ( y)$ and $\s^g (x,y) = 0$ otherwise. Here, the length $\ell(c) $ of a $C^1$ causal curve $c: [a,b] \rightarrow X$ is $\ell(c):= \int_{a}^b \sqrt{- g(c'(t) , c'(t))} dt $. If $X$ is a g.h. manifold-with-boundary, then $\s^g $ is finite and continuous (\cite{BEE}). Analogously to the Gromov compactification of locally compact separable length metric spaces, where, following an approach by Kuratowski, a point $x$ is identified with the distance $d(x, \cdot)$ (see e.g. \cite{FHS}, Sec.~4.1) we apply the same procedure to $ \sigma^g_x:= \sigma^g(x, \cdot) $. For $f \in M:= \{ f: \R \rightarrow \R \vert f {\rm  \ measurable \ and \  locally \  essentially \ bounded} \} $ and $p \in [1, \infty]$, we define pseudo-metrics $d_{g,f,p} (x,y) := \vert f \ci \sigma^g_x  - f \ci \sigma^g_y \vert_{L^p (X)} $. Compactness of each $J(x)$ implies $d_{g,f,p} < \infty$. We apply the usual functor $\la$ from metric spaces to length spaces (revised in Sec. \ref{Ansatzes}) to $d_{g,f,p}$. Two examples are:

\begin{itemize}

\item The {\em Noldus metric} $ d_{g, {\rm N}} := d_{g, \vert \cdot \vert, \infty}$, which never makes $X$ a length metric space (\cite{Noldus}, Th.6). In Sec. \ref{Ansatzes} we show that $(X, d_{g, {\rm N}}  )$ is in general not locally connected by rectifiable paths, thus $\la (d_{g, {\rm N}})$ does not generate the given topology: {\em it does not have sufficiently many geodesics}. 

\item The {\em Beem metric} $d_{g,{\rm B}}:=  d_{g, \chi_{]- \infty, 0[}, 1}+ d_{g, \chi_{]0, \infty[}, 1} $ ($\chi_A$ being the characteristic function of $A \subset \R$) satisfies $ d_{g,{\rm B}}(p,q)= \vol (J^+(p) \triangle J^+(q)) +  \vol (J^-(p) \triangle J^-(q))$ where $\triangle$ denotes the symmetric difference. The topology defined by this metric (originally suggested by Beem) is one of the two used in definitions of the future causal boundary (\cite{Beem}, \cite{oM-FC}). We show in Section \ref{Ansatzes} that this metric has splitting geodesics, so the same holds for $\lambda (d_{g,{\rm B}}) $, having the same geodesics, so {\em  for our purposes, $d_{g,{\rm B}} $ has too many geodesics.}

\end{itemize}

\V To define an interesting object in between those extrema, from now on we choose $p=2$, and let $h: \R \rightarrow \R$ with $h(x) := x^4 \ \forall x \in \R$, its crucial properties (see Sec. \ref{Proofs}) being $h \vert_{\R \setminus \{ 0 \}}, h'' \vert_{\R\setminus \{ 0\}} >0$ and $h(0) = h'(0) = h''(0) = h^{(3)} (0) = 0 $. Furthermore, for $r \in [-1,1]$ we define $f_r \in M$ by $f_r(x) := (\frac{1}{2}  + \frac{r}{2} \cdot \sgn (x) ) \vert x \vert \cdot x^3 \ \forall x \in \R$). We define $\Phi^{g,f}: x \mapsto f \ci \s_x^g$, $\Phi^g:= \Phi^{g,h}$ and $\Phi^g_r:= \Phi^{g, f_r}$, $d^g:= d_{g,h,2}$ and $d^g_r:= d_{g,f_r,2}$. For each $f$, the metric $\la ( d_{g,f,2} )  $ is the metric of the pullback of the scalar product of $L^2(X)$ via $\Phi^{g,f}$. In Sec. \ref{Proofs} we show:
 
\begin{Theorem}
\label{Main1}
Let $n \in \N $ and $X \in {\rm Obj} (C^-)$.
\begin{enumerate}[(i)]
\item $ \Phi^{g}: x \mapsto h \ci \sigma_x^g $ is a $C^2$ injective immersion of $X$ into $L^2(X)$, closed if $k_1 (X) , \Gamma (X) >0$.
\item $F: (X,g) \mapsto (X, \la (d^{g}_0) = (\Phi^{g} )^* (\langle \cdot ,  \cdot \rangle_{L^2(X)} ), (d^g_{-1/2}, d^g_0, d^g_{1/2})  )$ is an injective functor $C_n^- \rightarrow K_n$ whose push-down to isometry classes is injective, too. 
\item $  \forall \  b \in \R^7: b_1 \geq 0 \Rightarrow \exists \ a \in \R^4  \forall (X,g) \in C_n^-(b) : 
 (X, (\Phi^g)^* \langle \cdot , \cdot \rangle ) \in C_n^+ (a) $. 
 \item $ \forall \  b \in \R^7  : b_6 \leq b_4 \land b_3 \leq b_5 \Rightarrow C_n^-(b) $ contains an open set in $C_n^-$. 
\end{enumerate}
\end{Theorem} 

\V The last item shows that the range of applicability of the theorem is quite large. To prove the third item, we use certain links between second fundamental form, intrinsic and extrinsic diameter for submanifolds of Hilbert spaces that might be of independent interest (Theorems \ref{Hilbert-Ext-Int} and \ref{Hilbert-refined}). After establishing this link to Riemannian geometry, via a Riemannian result of Wong \cite{Wong} we obtain:

\begin{Theorem}
\label{Main2}
For every $n \in \N$ and for all $ b \in \R^7$ with $b_1 \geq 0$, there are only finitely many homeomorphism classes of compact Cauchy slabs in $C_n^-(b)$.
\end{Theorem}
 
\V Theorem \ref{Main2} is, to the author's knowledge, the first known finiteness result in Lorentzian geometry. A comparison with other approaches to this topic can be found in the last section.

\bigskip

 
{\bf Acknowledgment:} The author thanks Martin G\"unther for helpful comments on a first version. 
 
%


\section{Preliminaries on $\la$ and metrics for $p \neq 2$}
\label{Ansatzes}

We first revise well-known facts: Let $(Y,d) $ be a metric space, let $c \in C^0 ([a,b], Y)$. A {\em partition of $[a,b]$} is a finite subset $Z= \{ y_i \vert i \in \N_N \}$ of $[a,b]$ with $a,b \in Z$, which we always number monotonically. We denote the set of partitions of $[a,b]$ by $P(a,b)$ and put $\ell_Z (c) := \sum_{i=1}^N d(c_{i-1}, c_i)$ and $ (L (d))(c) := \sup \{ \ell_Z(c) \vert Z \in P(a,b) \}$, whose finiteness defines the set of {\em rectifiable curves}, containing all Lipschitz curves. Conversely, there is a map $K$ from length structures to $\infty$-metric spaces given by $ K(Y,\ell) := (Y,d) $ with $d(x,y) := \inf \{ \ell(c) \vert c: x \leadsto y\}$, then $ \lambda:= K \ci L$ is a non-injective functor from metric spaces to length spaces (with extended metrics, i.e. $\infty$-metric spaces), $\lambda \geq \1 $ (i.e., $\la (d) \geq d \ \forall d $), $ L \ci K = \1$ (see e.g. \cite{BBI}, Prop. 2.3.12), and $d$-geodesics are $\la(d)$-geodesics.

\begin{Theorem}[Noldus metric]
In $ (X:= [-1,1] \times \R, g= -dt^2 + ds^2)$ equipped with $d_{{\rm N}} := d_{g, \vert \cdot \vert , \infty}$, no continuous curve $c: [0,1] \rightarrow X$ for which $x_2 (c(0)) \neq x_2(c(1))$ is $d_{{\rm N}}$-rectifiable. 
\end{Theorem}

\begin{proof}
First assume w.l.o.g. that not both $c(0)$ and $c(1)$ are contained in the upper boundary $\{ 1 \} \times \R$ (otherwise reverse the time orientation). Then let, for $a >0$, 

$$A_a:= \{ t \in  [0,1] \vert x_1(c(t)) \leq 1-a \}  $$ 

and note that due to the continuity of $c$ there is a neighborhood $[0, t ]$ of $0$ contained in some $A_a$. With $u:= x_2 (c(0)), v:= x_2(c(t))$ we obtain $\ell(c\vert_{[0, t]}) \geq \ell(k)$, where $k: [0, t] \rightarrow X $ is defined by 

$$ k(s) := (1-a, (1-t^{-1}s) u + t^{-1} s v)  .$$

So due to rotational invariance everything boils down to calculating the $L^\infty$-distance $D(t)$ between $\s_{(1-a, 0)} $ and $\s_{(1-a, t)} $ which is $\sqrt{t(2a-t)}$. Then the argument is completed by calculating 

\begin{equation*}
n \cdot D(t/n)= \sqrt{t} \sqrt{2an-t } \rightarrow_{n \rightarrow \infty} \infty.  \qedhere
\end{equation*}

\end{proof}

\V Thus $X$ with the Noldus metric $d_N$ is not locally connected by rectifiable paths, so $\la (d_N)$ does not generate the topology of $X$. The same holds for $X= [0,1 ] \times {\mathbb{S}}^1$, with essentially the same proof.

\begin{Theorem}[Beem metric]
Every $g$-causal curve in $X$ is a geodesic for the Beem metric $d_{{\rm B}}$. 
\end{Theorem}

\begin{proof} We first note 
	
	$$ d_{{\rm B}} (x,y)= \vol (J^- (x) \triangle J^-(y) ) + \vol (J^+(x) \triangle J^+(y) ).$$
	
	For $ p \leq q \leq r$ we get 
	
	\bigskip
	
	$ J^\pm(p) \triangle J^\pm(r) = (J^\pm(p) \triangle J^\pm(q) ) \dot\cup (J^\pm(q) \triangle J^\pm(r)) $, so $d_{{\rm B}} (p,r) = d_{{\rm B}} (p,q) + d_{{\rm B}} (q,r)$. \end{proof}

\bigskip

Thus all $g$-causal curves are geodesics for $d_{{\rm B}}$ and $\la(d_{{B}})$, so they split, unlike in Alexandrov spaces.

%
%
%
%
%

\section{Proofs of Theorems 1 and 2, discussion}
\label{Proofs}

\subsection{Proof of Theorem 1 (i)}

%
%

The maps $\Phi^{g,f} : X \rightarrow L^2 (X)$ are obviously injective for every $f \in M$ with $f^{-1} (0) = \{ 0 \}$: A left inverse of $\Phi^{g,f}$ is given by $k \mapsto p$ if $p$ is the unique point in $\cl (k^{-1} (] - \infty, 0[)) \cap \cl (k^{-1} (]0, \infty[)) $. Keep in mind that the causal cut locus is of vanishing measure (\cite{BEE}), and that $y \in {\rm Cut}(x) \cap J^+(x) \Leftrightarrow x \in {\rm Cut} (y) \cap J^-(x)$.

\begin{Theorem}
\label{MetricAndHessian}
Let $f \in C^4(\R, \R)$ with $f^{-1} (0) = \{ 0 \}$ and $f^{(j)} (0) = 0 \forall j \in \N_3$, let $(X,g)$ be a g.h. manifold-with-boundary, then $\Phi^{g,f} : X \rightarrow L^2(X), x \mapsto f \ci \sigma^{g}_{x}$ is a $C^2$ embedding. For $w \in T_xX$, $\vert w \vert := (- g(w,w) )^{1/2}$ and $w(y):= \exp_x^{-1} (y)$ we get

\bean
d_x \Phi^{g,f}  \cdot v =   
 \Big( y \mapsto \begin{cases}  f'( \vert w(y)  \vert) \cdot \vert w(y) \vert^{-1} \cdot g( v, w (y) )  &  {\rm for \ }  w(y) \  {\rm  future \ timelike }\\ 
 	 - f'( - \vert w(y)  \vert) \cdot \vert w(y) \vert^{-1} \cdot g( v, w (y) )  &  {\rm for \ }  w(y) \  {\rm  past \ timelike }\\
  0 &{\rm else }
\end{cases}
\  \ \Big)
\eean

 (well-defined and smooth on the complement of the causal cut locus ${\rm Cut} (x)$ of $x$). The Hessian ${\rm Hess}_x (\Phi^{g,f}) $ of $\Phi^{g,f}$ at $x$ is defined on $X \setminus {\rm Cut}(x)   $ by (for $v \in T_xX$, $I_x:= \{v \in T_xX: g(v,v) <0 \}$): 
 
 \bean
{\rm Hess}_x (\Phi^{g,f}) \cdot (v \x v) =   
 \Big( y \mapsto   \begin{cases} f'(\vert w \vert )  \cdot \vert w (y) \vert^{-1} g( v, K_{v,y}'(0) ) + f''(\vert w (y) \vert) \cdot   (\vert w(y) \vert)^{-2} \cdot g( v, w )^2, &  w (y) \in I_x \\ 
  0 &{\rm else }
\end{cases}
\  \Big)
\eean

\V where $K_{v,y}$ is the Jacobi vector field along the unique maximal geodesic $c_y$ from $x= c_y(0)$ to $y=\exp_x(w) = c_y(1)$ with $K_{v,y}(1) = 0$ and $K_{v,y}(0)=v$. The pull-back Riemannian metric is 

\bea
\label{PullBackMetric}
(\Phi^{g,f})^* (\langle \cdot, \cdot \rangle_{L^2(X)})_x (a,b)   = \int_{J(x) } \vert w (y) \vert^{-1} (f'(\vert w (y) \vert ))^2 \cdot  g( a , w(y) ) \cdot g(  b , w(y) )  dy 
\eea

\end{Theorem}

{\bf Remark.} Despite of being $C^2$, in general $\Phi^{g,f}$ does {\em not} take values in $W^{1,p} (X)$.
 
\begin{proof} As $g$ and $f$ are fixed, we write $\Phi = \Phi{g,f}$. Let $x \in X$. We will frequently need an inverse of $\exp_x$, which is a priori only defined on a subset $ \exp_x (S_x) \subset X$ of full Lebesgue measure, which we complement by $\exp^{-1}_x (q) := 0 \in T_xX$ for all $q \in X \setminus \exp_x (S_x)$. As $\exp_x^{-1} $ is smooth on $\exp_x (S_x)$, $\exp_x^{-1}$ is continuous a.e. and is therefore measurable. As it is furthermore bounded and $\exp_x^{-1} (X) \cap I_x$ is compact, all involved functions are square-integrable. By the assumptions on $f $, the same is true of the expressions given for $d_x \Phi^{g,f} $ and $ \Hess _x (\Phi^{g,f}) $: To see this, first assume that $y = \exp _p(w) \gg x$. Then the first variational formula implies the formula for $d \Phi$, correspondingly for $x \gg y$.

For the Hessian ${\rm Hess} \Phi (V,W) := V(W ( \Phi)) - (\nabla_V W ) (\Phi) $ of $\Phi^{g,f}$ taking values in a Hilbert space with its trivial connection, which is a tensorial symmetric bilinear map, by polarization we only have to compute its values on $V=W=e_i$ where $e_i$ are $p$-synchronous vector fields such that $e_i (p)$ is a $g$-orthonormal basis of $T_pX$.

The second variational formula (\cite{ONeill}, Prop. 10.8.) determines the  second derivative of the length $l (c)$ for a geodesic variation around a non-null geodesic $c$ of signature $\e$ and speed $\kappa$ in the direction $V$ (which is a Jacobi vector field along $c$) as

\bean
\label{TransverseAcceleration}
V(V(l(0))) =   \frac{\e}{\kappa}  g(\nabla^*_{\partial t} V^\perp (t), V^\perp (t) ) \vert_a^b +  g(c'(t) , A(t)) \vert_a^b 
\eean

where $\nabla^*$  denotes the pull-back connection and $A$ is the transverse acceleration $\nabla^*_{\partial s} \partial_s F(0, \cdot)$ for the variation $F: (- \e, \e)_s \times [a,b]_t \rightarrow X $. For a synchronous vector field we have $A=0$.

The contribution of the shifting of the boundary of $J(x)$, which is of Lebesgue measure $0$, vanishes, because $\frac{f''(s)}{s^2} , \frac{f'(s)}{s^3}  \rightarrow_{s \rightarrow 0} 0$ and due to the formulas $d \Phi (x) \cdot v = (f' \ci \s_x ) \cdot (\s_x ' \cdot v)$ (here the first variational formula implies that the second factor is proportional to $(\s_x)^{-1}$) and $ d^2 \Phi (x) \cdot (v \x v) = (f'' \ci \s_x) \cdot (\s_x' \cdot v)^{\x 2} + (f' \ci \s_x) \cdot (\s_x''(v \x v)) $ (here the second variational formula gives a pole of first order in $\s_x$ for the second factor in the second term, whereas by the first variational formula the second factor in the first term has a pole of second order in $\s_x$). This fact is opposed to the situation in \cite{CGM}, Lemma 3.1, e.g., where the boundary term is central.

%
%
%
%

Now linearity of the solution $K$ of the Jacobi equation in its values at the endpoints gives bilinearity in $v$ and we conclude $ \s^g \in C^2(X, L^2(X))$.

Injectivity of $d_x \Phi$ follows: Let $v \neq 0$, then we find $w \in J^+_x  \subset T_xX $ with $g_x(v,w) \neq 0 $ and $t \in ]0, \infty [$ with $t w \in U_x$, which is the open subset of $T_xX$ from which $\exp_x$ is a diffeomorphism, then $ d_x \Phi (v) : \exp_x (tw) \mapsto (\frac{1}{2} + \frac{r}{2})^2 \langle v, w \rangle \neq 0   $, that is, $(d_x \Phi) (v) \neq 0$.   

Finally, $\Phi$ is closed (and thus an embedding) for $\Gamma (X) >0, \csec  $ finite: First we see along the lines of Lemma \ref{Lemma-Gamma} that there is a global (on $X$) positive lower bound on $\vol (J(p))$; defining recursively $C_n \in \partial^- X$ with $ J^-(J^+(C_n)) \subsetneq C_{n+1}$ and $A_n := J^+(C_n)$ we get for any $p_n \in A_n, q_n \in A_{n+2} \setminus A_n$ that $J(p_n) \cap J(q_n) = \emptyset$, so that there is $U >0 $ s.t. the $d(p_n,q_n) \geq U \ \forall n \in \N$.  \end{proof}

\subsection{Proof of Theorem 1 (ii)}

We recall $d_r^g (p,q) := \vert f_r \ci \s^g_p - f_r \ci \s^g_q \vert_{L^2(X)}$ for $p,q \in X$. We now suppress the dependence of $g$ for a moment in our notation. The fact that $\la(d^g_0)$ is a Riemannian metric follows from the previous item. For $x \in X$, let $\s_x^+$ resp. $\s_x^-$ denote the positive resp. negative part of $\s_x$. With the one-parameter family 

$ f_r: x \mapsto ( \frac{1}{2} + \frac{r}{2} \cdot \sgn (x) ) \cdot x^3 \vert x \vert$,
 
 $d_r$ interpolates between the past metric (taking into account only the past cones) $d_{-1} $  with $d_{-1} (x,y) := \vert \vert  (\s_x^-)^4 -  (\s_y^-)^4 \vert_{L^2} $ for $f_{-1} = (1- \theta_0) \cdot \1 ^4 $ and the future metric $d_1$ (taking into account only the future cones)  for $f_1 = \theta_0 \1 ^4$, passing through $d_0 (x,y) = \frac{1}{2} \big\vert \sigma_x^3 \cdot \vert \sigma_x \vert -  \sigma_y^3 \cdot \vert \sigma_y \vert \big\vert_{L^2}$. 
 Whereas $d_{\pm 1}^g$ is a metric on $ X \setminus \partial^\pm X$, it vanishes identically on $\partial^\pm X \times \partial^\pm X$.  For fixed points $p$ and $q$ of $X$, let $V$ be the two-dimensional linear span ${\rm span} (u^+, u^-)$ of $u^+ := \s_p^+  - \s_q^+$ and $u^- := \s_p^- - \s_q^-$. Then the $L^2$ scalar product on $V$ is uniquely given by the corresponding quadratic form on three vectors any two of which are non-collinear, thus given by $d_{-1/2} (p,q), d_0(p,q)$ and $d_{1/2} (p,q)$. In other words, the datum of those three recovers the whole family $d_r$. Then we can identify future and past boundary as $\partial^\pm X = \{ x \in X \vert \exists y \in X \setminus \{ x \} : d_{\pm 1}^g(x,y) = 0 \}  $. In the following we want to recover the causal structure. For $p,q \in X$ we define $\la_p^\pm := (\s_p^\pm)^4$ and consider the following expressions  only depending on the metrics $d_r$: 

$$ (d_{\pm 1 } (p,q) )^2     = \langle \la_p^\pm - \la_q^\pm , \la_p^\pm - \la_q^\pm \rangle_{L^2}   $$
 
 $$ (d_0 (p,q))^2 = \langle \la_p^+ - \la_p^- - \la_q^+ + \la_q^- , \la_p^+ - \la_p^- - \la_q^+ + \la_q^- \rangle_{L^2} , $$

we  calculate

$$   (d_0 (p,q))^2 -  (d_{-1 } (p,q) )^2 - (d_1 (p,q))^2 = 2 \big( \langle \la_q^+ , \la_p^- \rangle_{L^2} + \langle \la_p^- , \la_q^+ \rangle_{L^2}  - \underbrace{\langle  \la_p^+, \la_p^- \rangle_{L^2}}_{=0}  -  \underbrace{\langle  \la_q^+, \la_q^- \rangle_{L^2}}_{=0}   \big) ,   $$

where the last two terms vanish due to causality of $X$. This implies

$$  (d_0 (p,q))^2 -  (d_{-1 } (p,q) )^2 - (d_1 (p,q))^2  \neq 0 \Leftrightarrow p \ll q \lor q \ll p .$$

as $ p \ll q \Leftrightarrow \langle \s_p^+ , \s_q^- \rangle_{L^2} \neq 0$. The distinction between the two relevant cases can be met by taking into account $p \ll q \Leftrightarrow \s_p^+ > \s_q^+ $, so we only have to pick a point $r \in \partial^+ X $ (for which $d_1(p,r ) = \vert \s_p^+ \vert_{L^2}$ and $d_1(q,r ) = \vert \s_q^+ \vert_{L^2}$) and then 

$$p \ll q \Leftrightarrow (d_0 (p,q))^2 -  (d_{-1 } (p,q) )^2 - (d_1 (p,q))^2  \neq 0 \land d_1(p,r) > d_1(q,r) . $$

The causal structure can be recovered from the chronological structure as usual by 

$$p \leq q \Leftrightarrow  \forall r \in X: (q \ll r \Rightarrow p \ll r) , $$

thus we can identify the future and past subsets. These in turn form a subbasis for the manifold topology and allow to identify the conformal structure (\cite{BEE}, p.6, Th.2.3, Cor.2.4, Prop. 3.11), so everything is reduced to reconstructing the volume form. 

\bigskip

Now let $g$ and $h$ be two Lorentzian metrics $g,h$ on $X$ with $F (X,g)= F(X,h)$. By the above we can conclude that $g$ and $h$ are conformally related to each other, by, say, $g= e^{2u} \cdot h$ for a smooth function $u$. Assume that $u$ does not vanish everywhere, then w.l.o.g. let $u (x)>0$ for some $x \in X$. Due to continuity, there is an open neighborhood $U$ of $x$ and a real number $\e >0$ such that $u (y) > \e >0 \ \forall y \in U$.  As $(X,g)$ is g.h., there are $p,q \in X$ with $J_g^+(p) \cap J_g^-(q) \subset U  $. Then 

$$  \int_{J_g^+ (p) \cap J_g^-(q)} (\s^g_p)^4 (x) \cdot (\s^g_q)^4  (x) d \vol^g (x) \geq (e^{\e})^{n/2} \cdot (e^\e)^8  \int _{J_h^+ (p) \cap J_h^-(q)} (\s_p^{h})^4 (x) \cdot (\s_q^{h})^4 (x) d \vol^h (x) , $$

\V  But  we can reconstruct $   \int_{J_g^+ (p) \cap J_g^-(q)} (\s_p^{g})^4 (x) \cdot (\s_q^{g})^4 (x) d \vol^g (x) $ from the given data as above and exclude $u \neq 0$. Thus $g=h$. Of course, $F$ is also well-defined and injective on isomorphism classes, as each morphism on the  right-hand side induces an isomorphism on the left-hand side. \hfill \qed

\subsection{Proof of Theorem 1 (iii)}

We first define $\Phi:= \pr_1 \ci F: (X, g) \mapsto (X, \la (d^g)) $. First of all, the first item of the theorem establishes that the maps $\Phi$ are closed embeddings, thus the pull-back metric is complete.

The next four lemmas show how $\Phi$ transfers uniform Lorentzian bounds to Riemannian bounds.  Recall that cospacelike sectional curvature is sectional curvature on planes with spacelike orthogonal complement, which in the Lorentzian case are the Lorentzian planes.

\begin{Lemma}
\label{LemmaSec}
$k_1 (\Phi (k_1^{-1} ( [0, \infty [ )))  \subset [0, \infty[$ 

(in other words: when $k_1(X,g) \geq 0$ then  $k_1 ((X, \la(d_0^g)))  \geq 0 $).

If $k_1(X) \geq 0 $ then the Hessian $H_X^L$ of $\Phi^{g,f} : X \rightarrow L^2(X) $ takes values in the cone $P(X)$ of positive functions on $X$. 
\end{Lemma}

\begin{proof} The Gauss equation for the immersion $\Phi_r^g: X \rightarrow  L:=L^2(X)$ with Hessian $H^L_X$ reads:

$$  0= \langle R^L (V,Y)Z, W \rangle = \langle R^X (V,Y) Z , W \rangle + \langle H_X^L(V,Z) , H_X^L(Y,W) \rangle - \langle H_X^L(Y,Z), H_X^L(V,W) \rangle     $$

Thus we have to find a bound of the Hessian given by Th. \ref{MetricAndHessian}, which is the integral over $g( v=K_v(0), K_v'(0) )$ for Jacobi fields with $K_v(1)=0$. 
Let $J$ be a Jacobi field along timelike geodesics $c: I \rightarrow X$, w.l.o.g. orthogonal to $c'$. The crucial term in the definition of the Hessian is $u(v,w):= \vert K_{vw} \vert '(0) $ where $K_{vw}$ is the Jacobi  field along $c: t \mapsto \exp_x (tw)$ with $K_{vw} (1) =0$, $K_{vw} (0) = v$. Now if $\csec (X,g) \leq -\mu \leq 0$ then with $\kappa:= \mu \cdot \vert w \vert $ we get $g(R(V, c') c' , V) \geq \kappa \langle V, V )  $, and 
 
\bea
\label{Rauch}
  g( J, J )'' + \kappa g( J, J ) = g( J'', J ) + \kappa g( J,J ) + 2 g( J', J' ) \geq  g (R(J, c') c', J ) + \kappa g( J,J ) \geq 0  . 
  \eea

 Thus, for $u:= g( J, J ) $ and $\mu = 0$ we get $u'' \geq 0$. Then Eq. \ref{Rauch} with $u(0) = g( v, v)   >0$ and $u(\vert w \vert ) = 0$ means $0 >u'(0) = g( J(0), J'(0) )$. So $H_X^L:= \Hess (\Phi^{g,f})$ takes values in $P$. The Gauss equation in $L^2(X)$ implies


$$ 0 = \langle R^X (V,Y) Z,W \rangle_{L^2} + \langle H_X^L (V,Z), H_X^L (Y,W) \rangle_{L^2} - \langle H_X^L (Y,Z), H_X^L (V,W) \rangle_{L^2}  , $$

As $\sec (A,B) = \frac{\langle R(A,B)B, A \rangle }{\langle A,A \rangle \langle B,B \rangle - \langle A,B \rangle ^2}$, so for an orthonormal pair $(A,B)$ we have 

\bean
\sec(A,B) = \langle R(A,B)B, A \rangle = \underbrace{\langle H_X^L (B,B), H_X^L (A,A)  \rangle}_{=: W \geq 0} - \underbrace{\langle H_X^L (A,B) , H_X^L (A,B) \rangle}_{\vert \cdot \vert \leq_{(*)} W} \geq 0
\eean

where the starred inequality is due to the fact that for an orthonormal basis $(e_i)_{i \in \N}$ of positive functions we have (defining $W_i (U,V) := \langle H_X^L (U,V), e_i \rangle_{L^2}   $):

\bean
\vert \langle H_X^L (A,B) , H_X^L (A,B) \rangle_{L^2} \vert &=& \vert \sum_i \langle H_X^L (A,B) , e_i \rangle_{L^2} \langle e_i , H_X^L (A,B) \rangle_{L^2} \vert \\
&\leq& \sum_i \vert W_i (A,B) \cdot W_i (A,B) \vert \\
&\leq& \sum_i \vert W_i (A,A) \vert \cdot \vert W_i (B,B) \vert \\
&=_{(**)}& \sum_i W_i (A,A) W_i(B,B) \\
&=& \langle H_X^L (A,A), H_X^L (B,B) \rangle_{L^2}  
\eean

where the equality (**) is due to positivity of the basis elements. \end{proof}

\bigskip


\begin{Lemma}
For every $a_2, a_4, a_5, a_7 \in \R$, $k_2$ is bounded above on $\Phi (C_n^- (\infty, a_2, \infty , a_4, a_5, \infty, a_7))$.
\end{Lemma}

\begin{proof} The difficulty here is that we do not have a two-sided bound on the sectional curvature (and thus the Hessian of $\Phi$). Otherwise, to bound the second fundamental form $S_{\partial X}^X$, we could use $ S_{\partial X}^X = S_{\partial X}^{L^2(X)} - S_X^{L^2(X)} $, the fact that the second fundamental form is the normal part of the Hessian, the triangle inequality and the fact that $\partial^\pm X $ are Cauchy surfaces whose second fundamental form w.r.t. $g$ is uniformly bounded, allowing us to bound the transverse acceleration in Eq. \ref{TransverseAcceleration}. So we have to find another way: We need an upper bound on the integral 

$$A := \int_X \underbrace{\vert w(y) \vert^k \langle w(y), \nu \rangle}_{=: f(y)} \cdot \underbrace{\langle K_{e_j, y} ' (0), K_{e_j, y} (0) \rangle}_{=: h(y)} dy. $$ 

We find such a bound not via pointwise bounds but using Stokes' Theorem (observe that the last factor in the integral changes sign) in Federer's version \cite[Th. 4.5.6. (5)]{Federer} for Lipschitz vector fields. In Federer's book the theorem is formulated for open subsets of Euclidean half-spaces, but we can directly transfer it to pseudo-Riemannian manifolds-with-boundary via standard techniques like partition of unity and transformation formulas for tensors, and in our application we easily see via the choice of a Cauchy temporal function that the involved vector fields {\em are} Lipschitz). Defining the vector field $V(y):= P_{c_y} (w_y) = c_y'(\vert w(y) \vert) = d_{w(y)} \exp \cdot w(y)$ (recall that $w$ is smooth a.e.), we see that the second factor $h$ in the above integral is the divergence of $V$:

$$  K_{e_j, y}' (0) = \frac{\nabla}{dt} K_{e_j, y} (0) = \frac{\nabla}{dt}  \frac{d}{ds} V(0,0) = \frac{\nabla}{ds}  \frac{d}{dt} V(0,0)  = \nabla_{e_j} (w(y)) ,$$

and $K_{e_j,y} (0) = e_j$, thus for a parallelly extension of the orthonormal basis $(e_0 = w(y), ..., e_n)$ denoted by the same symbol we get 

$$ g( K_{e_j, y}' (0) , K_{e_j,y} (0) )  = g( \nabla_{e_j} V(y) , e_j (y) ) , $$

and summing up gives the divergence of $V$: In consequence, we get 

\bea
\label{Stokes}
 A= \int_X  V( \vert w \vert^k \cdot g( w , \nu  ) ) + \int_{\partial^+ X } f(y) \cdot g( V(y) , \nu (y) )  . 
 \eea

The first term in Eq. \ref{Stokes} is tangential to $c$, non-geometric and bounded by the the bounds on volumes of the causal cones and $\cdiam$, the second one by the bounds on $\cdiam$ and the volume of causal cones intersected with $\partial X$. \end{proof}

\begin{Lemma}
	\label{Lemma-Gamma}
	$\forall a_1, a_3, a_6 \in \R \ \exists b \in \R: k_4 (\Phi (C_n^-(a_1, \infty, a_3, \infty, \infty, a_6, \infty))) \subset [e^b , \infty [$.
\end{Lemma}

\begin{proof} The Ricci curvature bound implied by the sectional curvature bound ensures that we find $E>0$ s.t. for all $(X,g) \in A$ and all $y \in X$ we have $\exp_y^* \vol \geq E \cdot \vol_{\kappa}$ (for $\vol_{\kappa}=  \exp_{\tilde{x}}^* \vol_{M_\kappa}  $ in the Lorentzian model space $M_\kappa$ of constant sectional curvature $\kappa$)
within the domain of injectivity, via the Ricatti equation (cf \cite{Treude}).  For $x \in X$ and $v \in T_xX$, we denote by $J^+_x$ resp. $I^+_x$ the causal resp. timelike future cone in $T_xX$ and for $v \in I_x^+$ define $J^{+,v}_x := \{ u \in J_x^+ \vert g( u,v)/g(v,v)   \leq 1 \land g( u, u)  \geq \frac{1}{3} g ( v, v )  \}$. By definition, there is $E_2>0$ s.t. for all $(X,g) \in A$ and for every point $z$ in $({\rm injrad}^+_g)^{-1}( ]b, \infty [)$ there is $w \in T_zM$ timelike future with $g(w,w) < -E_2 $ and $\exp \vert {V_w} $ is a diffeomorphism onto its image for $ V_w := \{ v \in T_zX \vert v {\ \rm future , \ } 0 \leq v \leq w  \} $.  Now let us consider Eq. \ref{PullBackMetric}. Let $\vol(U) > \Gamma/2 $ and ${\rm injrad}^+_g (y) > \Gamma/2$ for all $y \in U$, let $x \in U$. We consider a $g$-pseudo-ONB $e_0,... , e_n$ at $x $ where $ae_0 = w$ for some $a >0$. For $k \in \N$, $k \leq n$, we consider the open subsets $W_k := \{ w \in T_xX \vert g(w,w) <0, w_0<1, w_k>1/2   \} $. They have $A$-uniformly large $\vol_\kappa$ in $T_xX$. Thus $(\la(d_o^g))_{kk} \geq  C_4 (g^+_V)_{kk}$, where $V=e_0$ and $g^+_V$ is the metric obtained by a Wick rotation around the normalized vector field $V$, which has the same volume form as $g$. Taking together the estimates, we get a bound $\la (d_0^g)> C_3 \cdot g^+_V $ at every point $ x \in X$, so we find constants $C_1, C_2>0$ s.t. all $g$ with timelike sectional curvature $\geq s $, $\ric \leq R$ and $\Gamma \geq e^u$ satisfy $\vol(X, \la(d^g_0)) \geq C_1 \vol (U, g_V^+) >C_2$.  \end{proof}



\bigskip

\V For the proof of Lemma \ref{Diam-bound} below, we will need some variant of the following result about the connection between intrinsic and extrinsic diameter, which is of independent interest:

\bigskip

\begin{Theorem}
	\label{Hilbert-Ext-Int}
	Let $H$ be a Hilbert space.
	
	Let $U(r,s):= \{ c \in  C^2([0,b], H ) : b \in \R, \vert c' \vert =1,   d(c(0) , c(t)) <r , \vert c''(t) \vert < s   \}$.
	\begin{enumerate}
		\item  If $ s < \frac{2 \sqrt{2}}{3r} =: \rho $ then $ \forall c \in U(r,s) : l(c) < \frac{3}{2} r$. 
		\item Let $p \in H$ and let $M$ be a closed submanifold of $B(p,r/2) $. If $\sqrt{\langle S_M^H , S_M^H \rangle + \langle S_{\partial M}^M, S_{\partial M}^M \rangle \cdot \chi_{\partial M}}  < \frac{2 \sqrt{2}}{3r}$, then the intrinsic diameter $\diam (M)$ satisfies $\diam (M) < \frac{3}{2} r$.  
	\end{enumerate}
\end{Theorem}

\V {\em Proof of Theorem \ref{Hilbert-Ext-Int}.} To prove the first item, assume w.l.o.g. that $c $ is parametrized by arc length. The mean value theorem asserts that for $\Delta (t):=  c'(t) - c'(0) $ we get $\vert \Delta (s) \vert < t \cdot \rho$. For $\a := \arccos (\langle  c'(t) , c'(0) \rangle ) $ and $\beta := \a/2$, we get by elementary trigonometry $\sin \beta = \Delta (t) /2$, thus 

\bea
\label{Absch1}
\langle c'(t), c'(0)  \rangle = \cos (2 \beta) &=& \cos^2 \beta - \sin^2 \beta \nonumber  \\ 
&=& 1 - 2 \sin^2 \beta = 1 - 2 \frac{\Delta^2 (t)}{4} = 1 - \frac{\Delta^2 (t)}{2} \geq 1 - \frac{1}{2} \rho ^2 t^2 
\eea

and consequently, by the Cauchy-Schwartz inequality, we get

\bean
\vert c(t) - c(0) \vert \geq \langle c(t) - c(0) , c'(0) \rangle &\geq& \int_0^t (1 - \frac{1}{2} \rho^2 \tau^2) d \tau \\
&=& t - \frac{1}{6} \rho^2 t^3 =: h(t)  
\eean

We are interested in the question when $h(t)$ becomes greater than $r$ and therefore calculate the maximum locus of $h$. The equation $0= h'(b_0) = 1 - \frac{1}{2} \rho^2 b_0^2 $ is solved by $b_0 = \frac{\sqrt{2}}{\rho}$. Inserted into $h$ this gives 

$$ h(b_0 ) = \frac{\sqrt{2}}{\rho} - \frac{2 \sqrt{2}}{6\rho} = \frac{2\sqrt{2}}{3\rho} =r , $$ 

and thus, if $c$ is defined on $[0, \sqrt{2}/\rho] = [0, \frac{3}{2}r ]$, then $\vert c(\sqrt{2}/\rho) - c(0) \vert > r$. 

For the second assertion, let $x, y \in M$. Closedness of $M$ and completeness of $H$ imply that there is a geodesic $c$ from $x$ to $y$. Then the extrinsic diameter of $c(I)$ is $\leq r$ by the triangle inequality. As $c$ is geodesic im $M$, we have $\nabla_t c' (t) \in (TM)^\perp$ and $\langle \nabla_t c'(t) , Y \rangle = \langle S_M^H (c'(t) , c'(t) ) , Y\rangle \ \forall t \in I $, so $\vert c'' (t) \vert = \vert S_M^H (c'(t), c'(t)) \vert \leq \vert S_M^H \vert $. 
Then we apply the first part of the theorem.  The bound on the second fundamental form of the boundary is needed as the connecting curve could have a part along the boundary where it is in general not a geodesic and has two normal parts, one in $M$ and one in $H$, whose sum is $c''$. \hfill \qed

\bigskip

Now, the direct application of the preceding lemma in our context would impose too severe restrictions, as two-sided bounds on cospacelike sectional curvature imply that sectional curvature is constant (\cite{Harris}, Prop. A.1). Instead, we pursue a different course: We first observe that by the uniform upper estimates $e^{a_4} $ resp. $ e^{a_3} $ on volume resp. timelike diameter on a subset $A \subset C$, we have that $\Phi^{g,f}(X) \in P \cap (v - P)$, for $v: X \rightarrow \R, x \mapsto e^{a_3}$. Furthermore, $L:= P \cap (v-P)$ is compact, as explained below. Then the following theorem ensures that the set $ \{ \diam ( \Phi^{g,f} (X), \la (d^{g,f})) \vert g \in A \}$ is bounded:

\begin{Theorem}
	\label{Hilbert-refined}
	Let $H$ be a separable Hilbert space and let $L \subset H$ be compact. Let $K$ be a convex self-dual cone in $H$ and $v \in K$. Then there is $D>0$ such that for all curves $c \in C^2( [a,b] ,L)$ parametrized by arc length with $ c' ([a,b]) \subset v - K$ and $c''([a,b]) \subset K$ we have $l (c) \leq D$. 
\end{Theorem} 

\V{\em Remark.} The hypothesis of self-duality is indispensable here, as otherwise for $v \in H \setminus \{ 0 \}$ and $K:= \{ u \in H \vert \langle u,v \rangle >0\}$ any periodic curve (or one contained in a finite-dimensional subspace) in the affine hyperplane $\{ u \in H \vert \langle u,v \rangle = \frac{1}{2} \langle v,v \rangle \} $ satisfies the hypotheses. Also $D$ has to depend on both $K$ and $L$, see the example given $H:= \R^2 \ni v$, $K:= (0, \infty [ ^2$, $L:= K \cap (v-K)$ and the curves $c$ being restrictions of $k: ]0, \infty [ \rightarrow L, k(t) := (t, 1/t) $.

\begin{proof} Let $n \in \N \cup \{ \infty \}$, $\N_{n} := \{ m \in  \N \vert m \leq n  \}$ and $\ell ^2 (n) := \{ a: \N_{n} \rightarrow \R \vert \sum_{k=1}^n a_k^2 < \infty  \}$. Then for any self-dual cone in the $n$-dimensional separable Hilbert space there is a unitary map $U: H \rightarrow \ell^2 (n)$ with $U(K) := \{ a \in \ell^2 (n) \vert a(k) \geq  0 \ \forall n \in \N_n \} \subset [0, \infty]^n $. Thus $U(L)$ is compact, in particular contained in $B(0,r)$ for some $r >0$, and there is an orthogonal projection onto $U(K)$ by just taking the positive part $w_+ \in U(K)$ of any vector $w$, and with $w_- := w_+ - w \in U(K)$ we get $w= w_+ - w_-$ and $\langle w_+, w_- \rangle = 0$ (in the example of $H= L^2 (X)$ this coincides with taking pointwise the positive and negative part of the functions). For all $u \in v-K $ there is some $E \in ]0, \infty[ $ with $E \cdot \vert \vert u \vert \vert \leq \langle u, v \rangle \leq E \cdot \vert \vert u \vert \vert $. Thus, by arc-length parametrization of $c$ and orthogonality of the projection on $K$, we see

\bea
1 = \langle c', c' \rangle = \langle c_+' , c_+' \rangle+ \langle c_-', c_-'\rangle \leq E^{-2} (\langle c_+', v \rangle^2 + \langle c_-', v \rangle^2) ,
\eea     

which implies $\forall t \in [a,b]:  (\langle c_+'(t) , v \rangle \geq E/\sqrt{2} ) \lor  (\langle c_+'(t) , v \rangle \geq E/\sqrt{2} ).$

With $L:= l(c) = b-a$, Cauchy's mean value theorem implies existence of some $t \in ]a,b[ $ with 

\bea
\label{BothParts}
l:= l(c) = \langle c_+'(t) , v \rangle - \langle c_-'(t) , v \rangle = \langle c'(t) , v \rangle  = \frac{1}{l} ( \langle c(b) , v  \rangle - \langle c(a), v \rangle  ) \leq 2 r \vert \vert v \vert \vert/l.
\eea

Now choose $D_0 \in \R$ such that

\bea
\forall l \geq D_0: \frac{E}{\sqrt{2}} - \frac{2}{l} r \vert \vert v \vert \vert > E/2 .
\eea 

Then Eq. \ref{BothParts} implies that {\em both} parts of the derivative satisfy a lower estimate: $(\langle c_+' (t) , v \rangle \geq E/2) \land (\langle c_-' (t) , v \rangle \geq E/2)$.

Assume that $c$ has length $\geq l$, then one of the two subintervals $[a,t]$ and $[t, b]$ has length $\geq l/2$. Assume first that this is true for $[a,t]$. Then

\bean
\langle c'(t), c_-'(t)  \rangle = - \langle c_- '(t) , c_-'(t) \rangle \leq - E^{-1} \langle c_-'(t) , v \rangle 
\leq - E^{-1} \cdot E/2 = 1/2.
\eean 

and therefore 

$$ \vert \vert c(0) - c(t) \vert \vert \geq \langle c(t) , c_- (t) \rangle - \langle c(0) , c_- (t) \rangle = \int_0^t \langle c'(s) , c_- (t) \rangle ds \geq \frac{l}{2} \cdot \frac{1}{2} $$

(here the last inquality holds as for $ u(s):= \langle c(s), c_-' (t) \rangle$ we get $u'(s) \geq 1/2$ and $c'' (s) \in K$ implies $u''(s) <0$ and thus  $ \langle c'(s) , c_-' (t) \rangle \leq \langle c'(t) , c_-' (t) \rangle \leq 1/2  $). Let $U:= \max \{ \vert \langle w, c_-' (t) \rangle \vert : w \in L \}  \in [0, \infty [ $, then $D:= \max \{ D_0, 8U \}$ implies, by convexity of $\langle \cdot, c_-'(t) \rangle $ along $c \vert_{[0,t]}$, that for all $l >D$ we have $\vert \vert c(0) - c(t)  \vert \vert \geq \frac{1}{2} \cdot \frac{l}{2} >2U $, in contrast to the definition of $U$. If the other subinterval was of length $\geq l/2$, then in the last computation we replace $c'_- (t)$ with $c_+' (t) $. \end{proof}

\begin{Lemma}
	\label{Diam-bound}
$\forall a_2, a_3, a_4 , a_6 \in \R \exists b \in \R:  k_4(\Phi(C_n^- (0, a_2, a_3, a_4, , \infty, a_6, \infty))) \subset (0, b].$

\end{Lemma}

\begin{proof} The bounds on the volume and the timelike diameter of $(X,g)$ imply a bound on the extrinsic diameters. Then the condition $H_X^L \in P$ on the Hessians $H_X^L$ of $\Phi^{g,f}: X \rightarrow L:= L^2(X)$ found in Lemma \ref{LemmaSec} together with the condition $v - H_X^L (p)  \in P \ \forall p \in \partial^- X$ bound the intrinsic diameter in terms of the extrinsic one by means of Theorem \ref{Hilbert-refined}: Let $x,y \in \Phi^{g,f} (X)$, then by completeness of $X$ there is a geodesic $c$ (modulo $\partial X$) from $x$ to $y$, whose second covariant derivative is the sum of a part normal to $\Phi^{g,f} (X)$, which is in the image of the Hessian taking values in $P$ and on the boundary of another part collinear to the exterior normal, (which is also a nonnegative function, as varying a point to the exterior increases $\s^+$), thus $\langle c'' (t) , w \rangle \leq 0  $. To show the last requirement $c'(s) \in v - K$, we first notice that arc-length parametrization of $c$ implies that for $y \gg c(s) $ and for any orthonormal basis $w(y) = e_0, e_1, ... , e_n$ at $c(s)$ we get as in the Lemma before that we find a uniform bound on $g(e_j, c'(s))$, and then we have to find a pointwise bound on $y \mapsto 4 \vert w(y) \vert^3 \cdot \langle v, w \rangle $ uniform in $C_n (0, a_2, a_3, a_4, a_5 , \infty) $. Indeed, for $w:= w_0 e_0 + W$ and $ v:= v_0 e_0 + W$ (where $V, W \perp e_0 $) we calculate $\vert w \vert^3 \cdot g( v,w ) = (\underbrace{w_0^2 - g(W,W)}_{>0} )^{3/2} \cdot (v_0 w_0 + \underbrace{g(V,W)}_{\leq \sqrt{g(V,V) \cdot g(W,W)}})\leq \vert w \vert^4 \cdot \underbrace{\vert v_0 + \sqrt{g(V,V)} \vert}_{< \sqrt{2} \vert v \vert_+} $, which yields the desired uniform upper bound.   \end{proof}

%

\subsection{Proof of Theorem 1 (iv)}
Let $S:= S_r:= \mathbb{S} (r) $ for $r > \pi^{-1}$, then $l:= \diam (S) >2$ and $\diam (S^{n-1}) = \sqrt{n-1} l  $. Let $X_r := ([-1, 1] \times S_r^{n-1} , g := -dt^2 + g^{n-1}_S)$, then $\cdiam (X_r)=2$, $\vol (X_r) = 2 l^{n-1}$, $\vol (\partial X_r) = l^{n-1}$, $\injrad^\pm (x) = \pm 1 \mp x_0  $, $V(b) := \vol ((\injrad^+)^{-1} ((b, \infty)) ) = \vol (X_r) - b \cdot \vol (X_r)/2$, and the claim is proven by perturbing around $X_r$ noticing that all involved geometric data are continuous real functions in the metric w.r.t. the $C^2$ topology.   

\subsection{Proof of Theorem 2}

%
%
%
%
%
%
%
%
%

The finiteness result of Theorem 2 now follows directly from Theorem \ref{Main1} and a result by Wong (where $C_n^a$ is denoted by $\mathcal{M} (n,K^-,\la^\pm, \vol\geq v>0,d)$, and for ease of comparison we note $K^- := a_1, \la^\pm := \pm e^{a_2}, v:= e^{a_3}, d:= e^{a_4}$):

\begin{Theorem}[\cite{Wong}, Th.1.4]
For every $n \in \N$ and every $a \in \R^4$, the number of homeomorphism classes in $ C^+_n(a)$ is finite. \hfill \qed 
\end{Theorem}

This concludes the Proof of Theorem \ref{Main2} \hfill \qed

%


\subsection{Discussion and indispensibility of the bounds in Th. \ref{Main1}}

First of all, it should be stressed that the considered Lorentzian bounds do not, by mere restriction to the Riemannian boundary metric, yield Riemannian bounds sufficient to apply Riemannian finiteness theorems. Furthermore, recall that in Cheeger-Gromov type finitenes results, the volume bound is needed to prevent dimension loss, as Perelman's theorem refers to Alexandrov spaces of Hausdorff dimension exactly $n$, whereas in Gromov's theorem all dimensions $\leq n$ are included. Due to dimensional homogeneity of Alexandrov spaces (\cite{BBI}, Th. 10.6.1) the dimension loss in a GH limit has to occur globally, if it occurs. Now, this dimension loss can actually happen in our context: For a sequence of thinner and thinner Lorentzian cylinders $([0,1] \times \mathbb{S}^1, -dt^2 \oplus \e  \cdot ds^2) $, the corresponding Riemannian cylinders become also infinitely thin in the limit $\e \rightarrow 0$, so there is dimensional degeneracy in this example. Therefore the lower bound on $\Gamma$ is indispensable in Theorems \ref{Main1} and \ref{Main2}. It could not be replaced with a mere bound on the total volume of $X$ either, which can be seen by considering a sequence of cylinders concentrating the volume at the boundary: There are functions $f_n: [-1, 1] \rightarrow [0, \infty [$ with $f_n (-x) = f_n(x) $ and $\lim_n f_n (\pm 1 ) = \infty$ such that $ \vol (X_n) \rightarrow_{n \rightarrow \infty}  1$ for $X_n := ([-1, 1] \times \mathbb{S}^1, -dt^2 + f_n ds^2)$ but still the limit of $ \pr_1 (F (X_n)) $ is one-dimensional. The same holds for flat cylinders, $f_n = n$ and $[0,1/n]$ instead of $[0,1]$. One way to obstruct this phenomenon is precisely the above combination $\Gamma$ of volume and injectivity radius.

 The bound on the volume is seen to be indispensable considering the examples of the semi-Riemannian products of long surfaces of genus $\rightarrow \infty$ and the Lorentzian line.

On one hand, it seems worthwhile to search for similar finiteness results using the notions of Lorentzian injectivity radius developed in \cite{Anderson}, \cite{LeFloch}, keeping in mind that it is fundamental to renounce the use of ad hoc temporal functions or, equivalently, ad hoc Riemannian metrics to obtain a genuinely Lorentzian result. On the other hand, it would of course be desirable to replace the condition on timelike sectional curvature with a Ricci bound $ c g (v,v) \leq \ric (v,v) \leq C g(v,v) $ on causal vectors $v$, more natural in the context of mathematical relativity. This looks feasible because of the definition of the metric as an integral over Jacobi fields, but here well-posedness of the initial-value problem for the vacuum Einstein equations puts strong limitations. In any case, replacing bounds on sec by bounds on Ric needs to be complemented including bounds on the injectivity radius (cf.\cite{Anderson}, Th.1)

One could be tempted to apply the {\em theory of optimal transport} between the signed length functions, but the Lorentzian version of optimal transport \cite{StefanSuhr2016} leads to a metric satisfying the {\em inverse} triangle inequality, disabling us to use Gromov compactness and Perelman's stability.

Trying to get a Lorentzian-Riemannian functor through a Wick rotation would require some natural temporal vector field, the only known such being gradients of CMC Cauchy temporal functions $t$, i.e. such that for all $a \in t(X) $ its $a$-level set is a Cauchy surface of constant mean curvature $a$. To construct those, e.g. as in \cite{Gerhardt2}, we need $\partial X$ to be CMC itself --- quite a nongeneric condition.

Finally, the fascinating framework of {\em Lorentzian length spaces} as in \cite{KunzingerSaemann} does not seem the right one for finiteness results, as its inverse triangle inequality again prevents us from getting Gromov precompactness --- but see the discussion in \cite{oM-Penrose}.



%
%
%


\begin{thebibliography}{99}


\bibitem{Anderson}
Michael T. Anderson: {\em Cheeger-Gromov Theory and Applications to General Relativity}, in: {The Einstein Equations and the Large Scale Behavior of Gravitational Fields}, ed. by Piotr T. Chru\'sciel and Helmut Friedrich, Birkh\"auser, Basel (2004).  arXiv:gr-qc/0208079


\bibitem{Beem}
John K. Beem: {\em A metric topology for causally continuous completions}, General Relativity and Gravitation 8, no 4, 245 --- 257 (1977)

\bibitem{BEE}
John K. Beem, Paul E. Ehrlich , Kevin L. Easley: {\em Global Lorentzian Geometry}, 2nd ed., Taylor \& Francis Inc. (1996)

\bibitem{mB}
Marcel Berger:{\em A Panoramic View of Riemannian Geometry}. Springer-Verlag (2002)

\bibitem{BBI}
Dmitri Burago, Yuri Burago, and Sergei Ivanov: {\em A Course in Metric Geometry}, 
American Mathematical Society (2001)

\bibitem{CGM}
Piotr T. Chruściel, James D. E. Grant and Ettore Minguzzi: {\em On differentiability of volume time functions}, Annales Henri Poincar\'e Vol. 17, Issue 10, 2801---2824 (2016). arxiv:1301.2909

\bibitem{Edwards}
David A. Edwards: {\em The Structure of Superspace}, in {\em Studies in Topology}, Academic Press (1975)

\bibitem{Federer}
Herbert Federer: {\em Geometric Measure Theory}. Springer-Verlag (1996)

\bibitem{FHS}
Jos\'e Luis Flores, Jonat\'an Herrera, Miguel S\'anchez: {\em Gromov, Cauchy and causal boundaries for Riemannian, Finslerian and Lorentzian manifolds}, Mem. AMS 226, No. 1064 (2013). arXiv:1011.1154


\bibitem{Gerhardt2}
Claus Gerhardt: {\em On the CMC foliation of future ends of a spacetime}. Pacific J. Math. 226, 297-308 (2006). arXiv:math/0408197

\bibitem{GPW}
Karsten Grove, Peter Petersen, and Jyh-Yang Wu: {\em Geometric finiteness theorems via controlled topology.} Invent. Math., 99:205 --- 213 (1991)

\bibitem{mG}
Mikhail Gromov: {\em Metric structures for Riemannian and non-Riemannian Spaces.} Progress in Mathematics vol. 152. Translated by Bates, Sean Michael. With appendices by M.Katz, P.Pansu, and S.Semmes 

\bibitem{Harris}
Steven G. Harris:  {\em A triangle comparison theorem for Lorentz manifolds}. Indiana Univ. Math. J. 31, no. 3, 289–308 (1982)


\bibitem{Kapovitch}
Vitali Kapovitch: {\em Perelman's stability theorem}, Surveys in Differential Geometry Vol. 11 (2006). arXiv:math/0703002

\bibitem{KunzingerSaemann}
Michael Kunzinger, Clemens S\"amann: {\em Lorentzian length spaces}, Ann. Global Anal. Geom. 54, no. 3, 399---447 (2018). arXiv:1711.08990

\bibitem{LeFloch}
Bing-Long Chen, Philippe G. LeFloch: {\em Injectivity Radius of Lorentzian Manifolds}, Commun. Math. Phys. 278, 679–713 (2008)




\bibitem{oM-FC}
Olaf M\"uller: {\em Topologies on the future completion}, arXiv:1909.03797

\bibitem{oM-Penrose}
Olaf M\"uller: {\em Functors in Lorentzian geometry --- three variations on  a theme}. arXiv:2205.01617

\bibitem{oM:CG}
Olaf M\"uller: {\em Cheeger-Gromov compactness for manifolds with boundary}.
arXiv.org: 1808.06458


\bibitem{Invitation}
Olaf M\"uller, Miguel S\'anchez: {\em An invitation to Lorentzian geometry}, Jahresbericht der Deutschen Mathematiker-Vereinigung, vol. 115, Issue 3-4, 153---183 (2014)

\bibitem{Noldus}
Johan Noldus: {\em The limit space of a Cauchy sequence of globally hyperbolic spacetimes}, Class.Quant.Grav. 21 (2004) 851---874. arxiv:gr-qc/0308075 

\bibitem{ONeill}
Barrett O'Neill: {\em Semi-Riemannian Geometry With Applications to Relativity}, Academic Press (1983)

\bibitem{StefanSuhr2016}
Stefan Suhr: {\em Theory of optimal transport for Lorentzian cost functions}, Mathematisches Institut (Universität Münster) (2018). arXiv:1601.04532

\bibitem{Treude}
Jan-Hendrik Treude, James D.E. Grant: {\em Volume Comparison for Hypersurfaces in Lorentzian Manifolds and Singularity Theorems}, Annals of Global Analysis and Geometry 43 (2013), 233---251. arXiv:1201.4249

\bibitem{aT}
Anton Tuzhilin: {\em Who invented the Gromov-Hausdorff distance?}, arXiv:1612.00728

\bibitem{Wong}
Jeremy Wong: {\em An extension procedure for manifolds with boundary}, Pacific Journal 235, 173 (2008)


\end{thebibliography}
\end{document}